\documentclass[letterpaper,11pt,oneside,reqno]{amsart}
\usepackage{bm}
\usepackage{mathrsfs}
\usepackage{yfonts}
\usepackage{amsfonts,amsmath, amssymb,amsthm,amscd}
\usepackage[height=9.6in,width=5.95in]{geometry}
\usepackage[mpexclude,DIV13]{typearea}
\usepackage{verbatim}
\usepackage{hyperref}
\usepackage{graphicx}
\usepackage[latin1]{inputenc}
\usepackage{latexsym}
\usepackage{lscape}
\usepackage{epsfig}
\usepackage{subfigure}
\include{bibtex}

\usepackage{setspace}

\begin{document}
\bibliographystyle{alpha}

\newcommand{\e}[0]{\epsilon}
\newcommand{\I}{{\rm i}}
\newcommand{\EE}{\ensuremath{\mathbb{E}}}
\newcommand{\PP}{\ensuremath{\mathbb{P}}}
\newcommand{\N}{\ensuremath{\mathbb{N}}}
\newcommand{\R}{\ensuremath{\mathbb{R}}}
\newcommand{\C}{\ensuremath{\mathbb{C}}}
\newcommand{\Z}{\ensuremath{\mathbb{Z}}}
\newcommand{\Q}{\ensuremath{\mathbb{Q}}}
\newcommand{\T}{\ensuremath{\mathbb{T}}}
\def \Ai {{\rm Ai}}
\newtheorem{theorem}{Theorem}[section]
\newtheorem{partialtheorem}{Partial Theorem}[section]
\newtheorem{conj}[theorem]{Conjecture}
\newtheorem{lemma}[theorem]{Lemma}
\newtheorem{proposition}[theorem]{Proposition}
\newtheorem{corollary}[theorem]{Corollary}
\newtheorem{claim}[theorem]{Claim}
\newtheorem{experiment}[theorem]{Experimental Result}

\def\todo#1{\marginpar{\raggedright\footnotesize #1}}
\def\change#1{{\color{green}\todo{change}#1}}
\def\note#1{\textup{\textsf{\color{blue}(#1)}}}

\theoremstyle{definition}
\newtheorem{rem}[theorem]{Remark}

\theoremstyle{definition}
\newtheorem{com}[theorem]{Comment}

\theoremstyle{definition}
\newtheorem{definition}[theorem]{Definition}

\theoremstyle{definition}
\newtheorem{definitions}[theorem]{Definitions}

\theoremstyle{definition}
\newtheorem{conjecture}[theorem]{Conjecture}
\title{Two ways to solve ASEP}

\author[I. Corwin]{Ivan Corwin}
\address{I. Corwin\\
  Massachusetts Institute of Technology, Department of Mathematics\\
  77 Massachusetts Avenue\\
  Cambridge, MA 02139-4307, USA}
\email{ivan.corwin@gmail.com}

\begin{abstract}
The purpose of this article is to describe the two approaches to compute exact formulas (which are amenable to asymptotic analysis) for the probability distribution of the current of particles past a given site in the asymmetric simple exclusion process (ASEP) with step initial data. The first approach is via a variant of the coordinate Bethe ansatz and was developed in work of Tracy and Widom in 2008-2009, while the second approach is via a rigorous version of the replica trick and was developed in work of Borodin, Sasamoto and the author in 2012.
\end{abstract}

\maketitle

\section{Introduction}

Exact formulas in probabilistic systems are exceedingly important, and when a new one is discovered, it is worth paying attention. This is a lesson that I first learned in relation to the work of Tracy and Widom on the asymmetric simple exclusion process (ASEP) and through my subsequent work on the Kardar-Parisi-Zhang (KPZ) equation. New formulas can enable asymptotic analysis and uncover novel (and universal) limit laws. Comparing new formulas to those already known can help lead to the realization that certain structures or connections exist between disparate areas of study (or at least can suggest such a possibility and provide a guidepost).

The purpose of this article is to describe the synthesis of exact formulas for ASEP. There are presently two approaches to compute the current distribution for ASEP on $\Z$ with step initial condition. The first (called here the {\it coordinate approach}) is due to Tracy and Widom \cite{TW1,TW2,TW3} in a series of three papers from 2008-2009, while the second (called here the {\it duality approach}) is due to Borodin, Sasamoto and the author \cite{BCS} in 2012.

The duality approach is parallel to an approach (also developed in \cite{BCS}) to study current distribution for another particle system, called $q$-TASEP. Via a limit transition, the duality approach becomes the replica trick for directed polymers. In fact, ASEP and $q$-TASEP should be considered as integrable discrete regularizations of the directed polymer model in which the replica trick (famous for being non-rigorous) becomes mathematically rigorous. Underlying the solvability of $q$-TASEP and directed polymers is an integrable structure recently discovered by Borodin and the author \cite{BorCor} called Macdonald processes (which in turn is based on the integrable system surrounding Macdonald symmetric polynomials). It is not presently understood where ASEP could fit into this structure, but the fact that the duality approach applies in parallel for ASEP and $q$-TASEP compels one to look for a higher structure which encompasses both.

\section{Current distribution for ASEP}

ASEP is an interacting particle system introduced by Spitzer \cite{Spitzer} in 1970 (though arising earlier in biology in the work of MacDonald, Gibbs and Pipkin \cite{MGP} in 1968). Since then it has become a central object of study in interacting particle systems and non-equilibrium statistical mechanics. Each site of the lattice $\Z$ may be inhabited by at most one particle. Each particle attempts to jump left at rate $q$ and right at rate $p$ ($p+q=1$), except that jumps which would violate the ''one particle per site rule'' are suppressed. We will assume $q>p$ and for later use call $q-p=\gamma$ and $p/q=\tau$ (note that $\gamma>0$ and $\tau<1$).

There are two ways of constructing ASEP as a Markov process. The ``occupation process'' keeps track of whether each site in $\Z$ is occupied or unoccupied. The state space is $Y=\{0,1\}^{\Z}$ and for a state $\eta=\{\eta_x\}_{x\in \Z} \in Y$, $\eta_x=1$ if there is a particle at $x$ and 0 otherwise. This Markov process is denoted $\eta(t)$.

The ``coordinate process'' keeps track of the location of each particle. Assume there are only $k$ particles in the system, then the state space $X_k=\{x_1<\cdots<x_k\}\subset \Z^k$ and for a state $\vec{x}=\{x_1<\ldots<x_k\}\in X_k$ the value of $x_j$ is the location of particle $j$. We call $X_k$ a Weyl chamber. Because particles cannot hop over each other, the ASEP dynamics preserve particle ordering. This Markov process is denoted $\vec{x}(t)$.

In this article, we will be concerned with the ``step'' initial condition for ASEP in which every positive integer site is initially occupied and every other site is initially unoccupied. In terms of the occupation process this corresponds to having $\eta_x(0) = \mathbf{1}_{x>0}$ (here and throughout $\mathbf{1}_E$ is the indicator function for event $E$). Let $N_x(\eta) = \sum_{y\leq x} \eta_y$ and note that $N_0(\eta(t))$ records the number of particles of ASEP which, at time $t$ are to the left of, or at the origin -- that is to say, it is the net current of particles to pass the bond $0$ and $1$ in time $t$.

\begin{theorem}\label{ASEPKPZ}
For ASEP with step initial condition and $q>p$,
\begin{equation*}
\lim_{t\to \infty} \PP\left(\frac{N_0(t/\gamma) - t/4}{2^{-1/3}t^{1/3}} \geq -s\right) = F_{\textrm{GUE}}(s),
\end{equation*}
where $F_{\textrm{GUE}}(s)$ is the GUE Tracy-Widom distribution.
\end{theorem}

\begin{rem}
The distribution function $F_{\textrm{GUE}}(s)$ can be defined via a Fredholm determinant as
\begin{equation*}
F_{\textrm{GUE}}(s) = \det(I - K_{\Ai})_{L^2(s,\infty)}
\end{equation*}
where Airy kernel $K_{\Ai}$ acts on $L^2(s,\infty)$ with integral kernel
\begin{equation*}
K_{\Ai}(x,y) = \int_{0}^{\infty} \Ai(x+t)\Ai(y+t) dt.
\end{equation*}
\end{rem}

For $q=1$ and $p=0$ this result was proved in 1999 by Johansson \cite{KJ} and for general $q>p$ it was proved by Tracy and Widom \cite{TW1,TW2,TW3} in 2009, and then reproved via a new formula by Borodin, Sasamoto and the author \cite{BCS} in 2012. This result confirms that for all $q>p$ ASEP is in the Kardar-Parisi-Zhang universality class \cite{KPZ} (see also the review \cite{ICreview}).

In order to prove an asymptotic result (such as above), it is very useful to have a pre-asymptotic (finite $t$) formula to analyze. If the formula does not increase in complexity as $t$ goes to infinity, there is hope to compute its asymptotics. Presently, there are two approaches to computing manageable formulas for the distribution of $N_0(t)$.

\section{The coordinate approach}
In \cite{TW1}, Tracy and Widom start by considering the ASEP coordinate process $\vec{x}(t)$ with only $k$ particles. In 1997, Sch\"{u}tz \cite{SchutzMastereq} computed the transition probabilities (i.e., Green's function) for ASEP with $k=2$ particles. The first step in \cite{TW1} is a generalization to arbitrary $k$. Let $P_{\vec{y}}(\vec{x};t)$ represent the probability that in time $t$, a particle configuration $\vec{y}$ will transition to a second configuration $\vec{x}$. As long as $p\neq 0$, it was proved in \cite{TW1} that
\begin{equation}\label{transitionform}
P_{\vec{y}}(\vec{x};t) = \sum_{\sigma\in S_k} \int \cdots \int A_\sigma \prod_{i=1}^{k} \xi_{\sigma(j)}^{x_j-y_{\sigma(j)}-1} e^{\e(\xi_j)t} d\xi_j,
\end{equation}
where the contour of integration is a circle centered at zero with radius so small as to not contain any poles of $A_{\sigma}$. Here $\e(\xi) = p \xi^{-1} +q\xi -1$ and
\begin{equation*}
A_{\sigma} = \prod\left\{S_{\alpha\beta}:\{\alpha,\beta\} \textrm{ is an inversion in }\sigma\right\}, \qquad S_{\alpha\beta} = -\frac{p+q\xi_{\alpha}\xi_{\beta} - \xi_{\alpha}}{p+q\xi_{\alpha}\xi_{\beta} -\xi_{\beta}}.
\end{equation*}

This result is proved by showing that that $P_{\vec{y}}(\vec{x};t)$ solves the master equation for $k$-particle ASEP
\begin{equation*}
\frac{d}{dt} u(\vec{x};t) = \big((L^{k})^* u)(\vec{x};t), \qquad u(\vec{x};0) = \mathbf{1}_{\vec{x}=\vec{y}}.
\end{equation*}
Here $(L^k)^*$ is the adjoint of the generator of the $k$-particle ASEP coordinate process (this just means that the role of $p$ and $q$ are switched in going between $L^k$ and $(L^k)^*$). For $k=1$, $L^1$ and $(L^1)^*$ act on function $f:\Z\to \R$ as
\begin{eqnarray*}
\big(L^1f\big )(x) &=& q\left[ f(x-1)-f(x)\right] + p\left[f(x+1)-f(x)\right],\\
\big((L^1)^*f\big)(x) &=& p\left[ f(x-1)-f(x)\right] + q\left[f(x+1)-f(x)\right].
\end{eqnarray*}
For $k>1$, the generator $L^k$ and its adjoint depend on the location of $\vec{x}$ in the Weyl chamber, reflecting the fact that certain particle jumps are not allowed near the boundary of the Weyl chamber.

Quoting a footnote in \cite{TW1}:
\begin{quotation}
The idea in Bethe Ansatz (see, e.g. \cite{LL,Sutherland,YangYang}), applied to one-dimensional $k$-particle quantum mechanical problems, is to represent the wave function as a linear combination of free particle eigenstates and to incorporate the effect of the potential as a set of $k-1$ boundary conditions. The remarkable feature of models amendable to Bethe Ansatz is that the boundary conditions for $k\geq 3$ introduce no more new conditions... The application of Bethe Ansatz to the evolution equation (master equation) describing ASEP begins with Gwa and Spohn \cite{GwaSpohn} with subsequential developments by Sch\"{u}tz \cite{SchutzMastereq}.
\end{quotation}

To see this in practice, assume that one wants to solve
\begin{equation*}
\frac{d}{dt} u(\vec{x};t) = \big((L^k)^* u\big)(\vec{x};t),\qquad u(\vec{x};0) = u_0(\vec{x})
\end{equation*}
for $\vec{x}$ in the Weyl chamber $X_k$.

\begin{proposition}\label{Prop1}
If $v:\Z^k\times \R_+\to \R$ solves the ``free evolution equation with boundary condition'':
\begin{enumerate}
\item For all $\vec{x}\in \Z^k$
$$\frac{d}{dt} v(\vec{x};t) = \sum_{j=1}^{k} \big([L^1]^*_j v\big)(\vec{x};t);$$
\item For all $\vec{x}\in \Z^k$ such that $x_{j+1}=x_{j}+1$ for some $1\leq j\leq k-1$,
$$p v(x_1,\ldots, x_j,x_{j+1}-1,\ldots, x_k;t) + q v (x_1,\ldots, x_j +1,x_{j+1},\ldots, x_k;t) - v(\vec{x};t)=0;$$
\item For all $\vec{x}\in X_k$, $v(\vec{x};0) = u_0(\vec{x})$;
\end{enumerate}
Then, for all $t\geq 0$ and $\vec{x}\in X_k$, $u(\vec{x};t) =v(\vec{x};t)$.
\end{proposition}
In (1) above, $[L^1]^*_{j}$ means to apply $(L^1)^*$ in the $x_j$ variable. In fact, some growth conditions must be imposed to ensure that $u$ and $v$ match (see Propositions 4.9 and 4.10 of \cite{BCS}) but we will not dwell on this presently.

This reformulation of the master equation involves only $k-1$ boundary conditions and is amendable to Bethe Ansatz -- hence one is led to postulate equation (\ref{transitionform}). It remains to check the ansatz (i.e., that $P_{\vec{y}}(\vec{x};t)$ solves the reformulated equation). The $A_{\sigma}$ are just right to enforce the boundary condition. The only challenge (which requires an involved residue calculation) is to check the initial data, since there are a total of $k!$ integrals.

The transition probabilities for $k$-particle ASEP is only the first step towards Theorem \ref{ASEPKPZ}. The next step is to integrate out the locations of all but one particle, so as to compute the transition probability for a given particle $x_m$. The formula for the location of the $m^{th}$ particle at time $t$ involves a summation (indexed by certain subsets of $\{1,\ldots, k\}$) of contour integrals. These formulas are a result of significant residue calculations and combinatorics.

At this point we are only considering $k$ particles, whereas for the asymptotic problem we want to consider step initial conditions. This is achieved by taking $y_j=j$ for $1\leq j\leq k$ and taking $k$ to infinity. After further manipulations, the $m^{th}$ particle location distribution formula has a clear limit as $k$ goes to infinity. This is the first formula for step initial condition and it is given by an infinite series of contour integrals.

In \cite{TW2} this infinite series is recognized as equal to a transform of a Fredholm determinant. By the simple relationship between the location of the $m^{th}$ particle of ASEP and $N_0(t)$ (defined earlier), this shows that
\begin{equation}\label{TWN0f}
\PP(N_0(t)=m) = \frac{-\tau^{m}}{2\pi \I} \int \frac{\det(I-\zeta K_1)}{(\zeta;\tau)_{m+1}} d\zeta,
\end{equation}
where the integral in $\zeta$ is over a contour enclosing $\zeta=q^{-k}$ for $0\leq k\leq m-1$ and $(a;\tau)_{n}=(1-a)(1-\tau a)\cdots (1-\tau^{n-1}a)$. Here $\det(I-\zeta K_1)$ is the Fredholm determinant with the kernel of $K:L^2(C_R)\to L^2(C_R)$ given by
\begin{equation*}
K_1(\xi,\xi') = q\frac{ e^{\e(\xi)t}}{p+q \xi\xi'-\xi},
\end{equation*}
and the contour $C_R$ a sufficiently large circle centered at zero.

There remains, however, a significant challenge to proving Theorem \ref{ASEPKPZ} from the above formula. As $m$ increases, the kernel $K_1$ has no clear limit, and the denominator term $(\zeta;\tau)_{m+1}$ behaves widely as $\zeta$ varies on its contour of integration. Much of \cite{TW3} is devoted to reworking the above formula into one for which asymptotics can be performed. This is done through significant functional analysis. The final formula, from which Theorem \ref{ASEPKPZ} is proved by asymptotics is that (leaving off the contours of integration)
\begin{equation}\label{Jeq}
\PP(N_0(t)\geq m) = \int \frac{d\mu}{\mu} (\mu;\tau)_{\infty}\det(I+\mu J),
\end{equation}
where the kernel of $J$ is given by
\begin{eqnarray*}
J(\eta,\eta')&=&\int \exp\{\Psi_{t,m,x}(\zeta)-\Psi_{t,m,x}(\eta')\}\frac{f(\mu,\zeta/\eta')}{\eta'(\zeta-\eta)}d\zeta,\\
f(\mu,z)&=&\sum_{k=-\infty}^{\infty} \frac{\tau^k}{1-\tau^k\mu}z^k,\\
\Psi_{t,m,x}(\zeta) &=& \Lambda_{t,m,x}(\zeta)-\Lambda_{t,m,x}(\xi),\\
\Lambda_{t,m,x}(\zeta) &=& -x\log(1-\zeta) + \frac{t\zeta}{1-\zeta}+m\log\zeta.
\end{eqnarray*}

\section{The duality approach}

Duality is a powerful tool in the study of Markov processes. It reveals hidden structures and symmetries of the process, as well as leads to non-trivial systems of ODEs which expectations of certain observables satisfy. In 1997, Sch\"{u}tz \cite{SchutzDuality} observed that ASEP is self-dual (in a sense which will be made clear below). The fact that duality gives a useful tool for computing the moments of ASEP was first noted by Imamura and Sasamoto \cite{IS} in 2011. In 2012, Borodin, Sasamoto and the author \cite{BCS} used this observation about duality, along with an ansatz for solving the duality ODEs (which was inspired by the work of Borodin and the author on Macdonald processes \cite{BorCor}) to derive two different formulas for the probability distribution of $N_0(t)$. The first was new and readily amendable to asymptotic analysis necessary to prove Theorem \ref{ASEPKPZ}, while the second was equivalent to Tracy and Widom's formula (\ref{TWN0f}).

To define the general concept of duality, consider two Markov processes $\eta(t)$ with state space $Y$ and $\vec{x}(t)$ with state space $X$ (for the moment we think of these as arbitrary, though after the definition of duality, we will take these as before). Let $\EE^{\eta}$ and $\EE^{\vec{x}}$ represent the expectation of these two processes (respectively) started from $\eta(0)=\eta$ and $\vec{x}(0)=\vec{x}$. Then $\eta(t)$ and $\vec{x}(t)$ are dual with respect to a function $H:Y\times X\to \R$ if for all $\eta\in Y$, $\vec{x}\in X$ and $t\geq 0$,
\begin{equation*}
\EE^{\eta}\left[ H(\eta(t), \vec{x})\right] = \EE^{\vec{x}}\left[ H(\eta, \vec{x}(t))\right].
\end{equation*}
One immediate consequence of duality is that if we define $u_{\eta}(\vec{x};t)$ to be the expectations written above, then
\begin{equation*}
\frac{d}{dt} u_{\eta}(\vec{x};t) = L u_{\eta}(\vec{x};t)
\end{equation*}
where $L$ is the generator of $\vec{x}(t)$ and where the initial data is given by $u_{\eta}(\vec{x};0) = H(\eta,\vec{x})$.

Sch\"{u}tz \cite{SchutzDuality} observed that if $\eta(t)$ is the ASEP occupation process and $\vec{x}(t)$ is the $k$-particle ASEP coordinate process with $p$ and $q$ switched from the earlier definition, then these two Markov processes are dual with respect to
\begin{equation*}
H(\eta,\vec{x}) = \prod_{j=1}^{k} \tau^{N_{x_j-1}(\eta)} \eta_{x_j}.
\end{equation*}
The generator of the $p,q$ reversed particle process $\vec{x}(t)$ is equal to $(L^k)^*$, as discussed earlier. Sch\"{u}tz demonstrated this duality in terms of a spin-chain encoding of ASEP by using  a commutation relation along with the $U_q[SU(2)]$ symmetry of the chain. A direct proof can also be given in terms of the language of Markov processes \cite{BCS}. When $p=q$, $\tau=1$ and this duality reduces to the classical duality of correlation functions for the symmetric simple exclusion process (see \cite{Lig} Chapter 8, Theorem 1.1).

As before, we focus on step initial condition, so that $\eta_x = \mathbf{1}_{x\geq 1}$. Duality implies that $u_{\textrm{step}}(\vec{x};t):=\EE^{\eta}\left[ H(\eta(t), \vec{x})\right]$ solves
\begin{equation}\label{utaustep}
\frac{d}{dt} u_{\textrm{step}}(\vec{x};t) = L^k u_{\textrm{step}}(\vec{x};t), \qquad u_{\textrm{step}}(\vec{x};0) = \mathbf{1}_{x_1\geq 1} \prod_{i=1}^{k} \tau^{x_{i}-1}.
\end{equation}

The above system is solved by
\begin{equation}\label{usoln}
u_{\textrm{step}}(\vec{x};t) = \frac{\tau^{k(k-1)/2}}{(2\pi \I)^k} \int \cdots \int \prod_{1\leq A<B\leq k} \frac{z_A-z_B}{z_A-\tau z_B} \prod_{j=1}^{k} h_{x_j,t}(z_j) dz_j,
\end{equation}
where
\begin{equation*}
h_{x,t}(z) = e^{\e'(z) t} \left(\frac{1+z}{1+z/\tau}\right)^{x-1} \frac{1}{\tau +z}, \qquad \e'(z) = -\frac{z(p-q)^2}{(1+z)(p+qz)},
\end{equation*}
and where the contour of integration for each $z_j$ is a circle around $-\tau$, so small as to not contain $-1$. In order to see this, we use the reformulation of the system (\ref{utaustep}) in terms of the free evolution equation with boundary condition with ASEP given earlier in Proposition \ref{Prop1}. Condition (1) is trivially checked since for each $z$, $\tfrac{d}{dt} h_{x,t}(z)=L^1 h_{x,t}(z)$. Condition (3) is checked via a simple residue calculation. Condition (2) reveals the purpose of the $\tfrac{z_A-z_B}{z_A-\tau z_B}$ factor. Applying the boundary condition to the integrand above brings out a factor of $z_j-\tau z_j$. This cancels the corresponding term in the denominator and the resulting integral is simultaneous symmetry and antisymmetry in $z_j$ and $z_{j+1}$. Hence the integral must equal zero, which is the desired boundary condition (2).

The inspiration for this simple solution to the system of ODEs came from analogous formulas which solve free evolution equations with boundary condition for various versions of the delta Bose gas (see Section \ref{replicas} for a brief discussion). For the delta Bose gas and certain integrable discrete regularizations, the formulas arose directly from the structure of Macdonald processes \cite{BorCor}. ASEP does not fit into that structure, but the existence of similar formulas suggests the possibility of a yet higher structure.

A change of variables reveals some similarities to the integrand in (\ref{transitionform}). Letting
\begin{equation}\label{xichange}
\xi_j = \frac{1+z_j}{1+z_j/\tau}
\end{equation}
we have that
\begin{equation*}
\frac{z_A-z_B}{z_A-\tau z_B} = q\frac{\xi_A-\xi_B}{p+q\xi_A\xi_B-\xi_B},\qquad h_{x_j,t}(z_j) dz_j = e^{\e(\xi_j)t} \xi_j^{x_j-1} \frac{d\xi_j}{\tau-\xi_j}.
\end{equation*}

The system (\ref{utaustep}) could also be solved via Tracy and Widom's formula (see formula \ref{transitionform} earlier) for the Green's function for $(L^k)^*$ (as suggested in \cite{IS}) but the resulting formula would involve the sum of $k!$ $k$-fold contour integrals. Symmetrizing (\ref{usoln}) via combinatorial identities, and making the above change of variables, one does recover that formula. The reversal of this procedure is a rather unnatural anti-symmetrization, which explains why (\ref{usoln}) was not previously known.

A suitable summation of $H(\eta,\vec{x})$ over $\vec{x}$ gives $\tau^{kN_{x}(\eta)}$. Using this, and formula (\ref{usoln}), \cite{BCS} proves that for ASEP with step initial condition,
\begin{equation}\label{21}
\EE\left[\tau^{kN_0(t)}\right] = \frac{\tau^{k(k-1)/2}}{(2\pi \I)^k} \int \cdots \int \prod_{1\leq A<B\leq k} \frac{z_A-z_B}{z_A-\tau z_B} \prod_{j=1}^{k} e^{\e'(z_j)t}\frac{dz_j}{z_j},
\end{equation}
where $N_0(t) = N_0(\eta(t))$ and where the contour of integration for $z_j$ includes $0,-\tau$ but not $-1$ or $\tau$ times the contours for $z_{j+1}$ through $z_{k}$. This is to say, that the contours of integration respect a certain nesting structure.

At this point the utility of having a single $k$-fold nested contour integral formula for the moments of $\tau^{N_0(t)}$ becomes clear. There are two ways to deform the contours of integration in (\ref{21}) so as to all coincide with each other. The first involves expanding them all to be a circle containing $-\tau$ and 0, but not $-1$. There are many poles encountered in the course of this deformation and the residues can be indexed by a partition. This leads to
\begin{equation}\label{Etn}
\EE\left[\tau^{kN_0(t)}\right]= k_\tau! \sum_{\substack{\lambda\vdash k\\ \lambda=1^{m_1}2^{m_{2}}\cdots}} \frac{1}{m_1!m_2!\cdots} \, \frac{(1-\tau)^{k}}{(2\pi \iota)^{\ell(\lambda)}} \int \cdots \int \det\left[\frac{-1}{w_i \tau^{\lambda_i}-w_j}\right]_{i,j=1}^{\ell(\lambda)} \prod_{j=1}^{\ell(\lambda)} e^{t\sum_{i=0}^{\lambda_j-1} \e'(\tau^i w_j)} dw_j,
\end{equation}
where $k_{\tau}!=(\tau;\tau)_{k} (1-\tau)^{-k}$ is the $\tau$-deformed factorial, and $\lambda=(\lambda_1\geq \lambda_2\geq \cdots \geq 0)$ is a partition of $k$ (i.e. $\sum \lambda_i=k$) with $\ell(\lambda)$ nonzero parts, and multiplicity $m_j$ of the value $j$. The structure of these residues is very similar to the string states indexing the eigenfunctions of the attractive delta Bose gas (see Section \ref{replicas}).

The final step in the duality approach is to use these moment formulas to recover the distribution of $N_0(t)$. This is done via the $\tau$-deformed Laplace transform Hahn \cite{Hahn} introduced in 1949. The left-hand side of the below equation is the transform of $\tau^{N_0(t)}$ with spectral variable $\zeta$.
\begin{equation}\label{23}
\EE\left[\frac{1}{(\zeta \tau^{N_0(t)};\tau)_{\infty}}\right]  = \sum_{k=0}^{\infty} \frac{\zeta^k \EE [\tau^{k N_0(t)}]}{(\tau;\tau)_{k}}.
\end{equation}
The right-hand side above comes from the left-hand side by expanding the $\tau$-deformed exponential inside the expectation (using the $\tau$-deformed Binomial theorem) and then interchanging the summation over $k$ with the expectation. This interchange of summation and integration is justified here for $\zeta$ small enough because $|\tau^{kN_0(t)}|\leq 1$ deterministically (in contrast to (\ref{reptrick}) Section \ref{replicas}).

Substituting (\ref{Etn}) into the series on the right-hand side of (\ref{23}) one recognizes a Fredholm determinant. The kernel of the determinant can be rewritten using a Mellin-Barnes integral representation and the result is (leaving off the contours of integration)
\begin{equation}\label{mbdet}
\EE\left[\frac{1}{(\zeta \tau^{N_0(t)};\tau)_{\infty}}\right] = \det(I+K_\zeta),
\end{equation}
where the kernel of $K_{\zeta}$ is
\begin{equation*}
K_{\zeta}(w,w') = \frac{1}{2\pi \I} \int \frac{\pi}{\sin(-\pi s)} (-s)^{\zeta} \frac{g(w)}{g(\tau^s w)} \frac{ds}{w'-\tau^s w}, \qquad g(w) = e^{\gamma t \frac{\tau}{\tau+w}}.
\end{equation*}

The $\tau$-Laplace transform can easily be inverted to give the distribution of $N_0(t)$ and asymptotics of the above formula are readily performed (see Section 9 of \cite{BCS}) resulting in Theorem \ref{ASEPKPZ}.

There is a second choice for how to deform the nested contours in (\ref{Etn}) to all coincide. The terminal contour of this deformation is a small circle around $-\tau$, and again there are certain poles encountered during the deformation. The combinatorics of the residues here is simpler than in the first case, and one finds the following Fredholm determinant formula
\begin{equation}\label{cauchydet}
\EE\left[\frac{1}{(\zeta \tau^{N_0(t)};\tau)_{\infty}}\right] = \frac{\det(I- \zeta K_2)}{(\zeta;\tau)_{\infty}}
\end{equation}
where the kernel of $K_2$ is
\begin{equation*}
K_{2}(w,w') = \frac{e^{\e'(w)t}}{\tau w-w'}.
\end{equation*}
Performing the change of variables (\ref{xichange}) and inverting this $\tau$-Laplace transform, one recovers Tracy and Widom's formula (\ref{TWN0f}). As in Tracy and Widom's work, this formula is not yet suitable for asymptotics and must be manipulated significantly to get to the form of (\ref{Jeq}).

\section{Duality approach as a rigorous replica trick}\label{replicas}

Besides the Sch\"{u}tz duality, Borodin, Sasamoto and the author discovered that ASEP is also self dual with respect to
\begin{equation*}
H(\eta,\vec{x}) = \prod_{j=1}^{k} \tau^{N_{x_j}(\eta)}.
\end{equation*}
For $k=1$ this shows that $\EE[\tau^{N_x(\eta(t))}]$ solves the heat equation with generator $L^1$. In fact, this is essentially G\"{a}rtner's 1988 observation \cite{G} that $\tau^{N_x(\eta(t))}$ solves a certain discrete multiplicative stochastic heat equation. A multiplicative stochastic heat equation has a Feynman-Kac representation which shows that the solution can be interpreted as a partition function for a directed polymer in a disorder given by the noise of the stochastic heat equation.

In 1997 Bertini and Giacomin \cite{BG} showed that under a certain ``weakly asymmetric'' scaling $\tau^{N_x(\eta(t))}$ converges to the solution to the continuum multiplicative stochastic heat equation (SHE) with space-time white noise $\xi(x,t)$:
\begin{equation*}
\frac{d}{dt} Z(x,t) = \frac{1}{2} \frac{d^2}{dx^2} Z(x,t) + Z(x,t)\xi(x,t).
\end{equation*}
This convergence result did not include when $\eta(0)$ is step initial condition and was extended to that case by Amir, Quastel and the author \cite{ACQ}. The corresponding initial data for the SHE is $Z(x,0)=\delta_{x=0}$ where $\delta$ is the Dirac delta function. The logarithm of the solution to the SHE (formally) solves the Kardar-Parisi-Zhang equation
\begin{equation}
\frac{d}{dt} h(x,t) = \frac{1}{2} \frac{d^2}{dx^2} h(x,t) + \frac{1}{2} \left(\frac{d}{dx} h(x,t)\right)^2 + \xi(x,t).
\end{equation}
See \cite{ICreview} for more details on the Kardar-Parisi-Zhang equation.

Duality of ASEP translates into the fact that the moments of the SHE solve the attractive one-dimensional imaginary-time delta Bose gas (Lieb-Liniger model with delta interaction) \cite{K}. Define $\bar{Z}(\vec{x};t) = \EE\left[Z(x_1,t)\cdots Z(x_k,t)\right]$ for $Z$ with $\delta_{x=0}$ initial data. Then $\bar{Z}$ solves the system
\begin{equation}\label{dbosegas}
\frac{d}{dt} \bar{Z}(\vec{x};t) = H_1 \bar{Z}(\vec{x};t),\qquad \bar{Z}(\vec{x};0) = \prod_{j=1}^{k} \delta_{x_j=0},
\end{equation}
where $H_{\kappa}$ is the Lieb-Liniger Hamiltonian with delta interaction with strength $\kappa\in \R$:
\begin{equation*}
H_{\kappa}= \frac{1}{2} \sum_{j=1}^{k} \frac{d^2}{dx_j^2} + \kappa \sum_{i<j} \delta_{x_i=x_j}.
\end{equation*}

The Lieb-Liniger model with delta interaction was the second system solved by the Bethe Ansatz (over 30 years after Bethe \cite{Bethe} solved the spin-$1/2$ isotropic Heisenberg model). This was accomplished by Lieb and Liniger in 1963 for the repulsive system ($\kappa<0$). A year later, McGuire similarly solved the attractive system ($\kappa>0$). In their context, solving the system meant writing down eigenfunctions for $H_{\kappa}$. The structure of the eigenfunctions for the repulsive versus attractive cases are different. In the attractive case there are extra eigenfunctions which are called ``string states'' due to the strings of quasi-momenta with which they are indexed (or physically corresponding to bound states of particle clusters). Completeness of these eigenfunctions was not shown until later \cite{Dorlas,HeckOp,TWBose,Oxford,ProSpoComp}.

For the purposes of understanding the moments of the SHE it is not necessary to diagonalize $H_\kappa$, but rather just to solve the system (\ref{dbosegas}) for $\kappa=1$. Just as with ASEP, this system can be written as a ``free evolution equation with boundary condition''. The free evolution is just according to the $k$ variable Laplacian and the boundary condition is that for all $1\leq j\leq k-1$,
\begin{equation*}
\left(\frac{d}{dx_{j}} - \frac{d}{dx_{j+1}} - \kappa\right) v(\vec{x};t)\big\vert_{x_j\to x_{j+1}} = 0.
\end{equation*}

This system can be solved via an analogous formula to (\ref{usoln}): For $x_1\leq \ldots \leq x_k$ and $\kappa\in \R$,
\begin{equation}\label{deltanested}
\bar{Z}(\vec{x};t)  = \frac{1}{(2\pi \I )^k} \int \cdots \int \prod_{1\leq A<B\leq k} \frac{z_A-z_B}{z_A-z_B-\kappa} \prod_{j=1}^{k} \exp\left\{\frac{z_j^2}{2} t + x_jz_j\right\}dz_j
\end{equation}
where the contour of integration for $z_j$ is along $\alpha_j+ i \R$ for any $\alpha_1>\alpha_2 + \kappa >\cdots \alpha_k + (k-1)\kappa$. When $\kappa<0$ all the $\alpha_j$ can be chosen as $0$ and hence the integral occurs on $\I\R$, whereas for $\kappa>0$ the contours must be spaced horizontally. In the $\kappa>0$ case, the contours can be deformed to $\I\R$. The singularities and associated residues encountered have a very similar structure to those seen earlier in (\ref{Etn}) in the context of the first ASEP contour deformation. The disparity between residue combinatorics accounts for the difference in the structure of the eigenfunctions and the occurrence of string stated for $\kappa>0$. In fact, Heckeman and Opdam's 1997 proof of the completeness of the Bethe Ansatz relied on a formula equivalent to (\ref{deltanested}).

Given expressions for all of the moments of the SHE, one wants to recover the distribution of $Z(x,t)$. Since $Z(x,t)$ is nonnegative, its Laplace transform characterizes its distribution. N\"{a}ively one writes
\begin{equation}\label{reptrick}
\EE\left[e^{\zeta Z(t,x)}\right] = \sum_{k=0}^{\infty} \frac{ \zeta^k \EE[Z(t,x)^k]}{k!}.
\end{equation}
However, the right-hand side is known to make no mathematical sense and the interchange of expectation and summation is totally unjustifiable. The moments of the SHE grow like $e^{ck^3}$ and thus the right-hand side is extremely divergent. One can see that cutting off the summation also fails to remedy the situation in any way.

What should be clear now is that ASEP is an integrable discrete regularization of the SHE (or equivalently the KPZ equation) and the duality approach to solving it is a rigorous version of the replica trick for the SHE. By taking the weakly asymmetric limit of the $\tau$-deformed Laplace transform formulas described above, one finds a Fredholm determinant formula for $\EE\left[e^{\zeta Z(t,x)}\right]$. This can be done from either the new formula (\ref{mbdet}) in \cite{BCS} or Tracy and Widom's formula (\ref{Jeq}). It appears that (\ref{mbdet}) is very amendable to asymptotic analysis.

Using (\ref{Jeq}), the derivation of the Laplace transform of $Z(t,x)$ involves extremely careful asymptotic analysis which was performed in 2010 rigorously by Amir, Quastel and the author \cite{ACQ} and independently and in parallel (though non-rigorously) by Sasamoto and Spohn \cite{SS1}. Very soon afterwards, Calabrese, Le Doussal and Rosso, as well as Dotsenko showed how to formally recover this Fredholm determinant from summing the divergent series on the right-hand side of (\ref{reptrick}). The formal manipulations of divergent series that goes into this can be see as shadows of the rigorous duality approach explained above for ASEP. It can also be seen as a shadow of a parallel duality approach for q-TASEP \cite{BorCor,BCS}, another integrable discrete regularization of the SHE.

\subsection{Acknowledgements}
The author was partially supported by the NSF through grant DMS-1208998, PIRE grant OISE-07-30136 as well as by Microsoft Research through the Schramm Memorial Fellowship, and by the Clay Mathematics Institute.


\begin{thebibliography}{alpha}
\bibitem{ACQ}
G.~Amir, I.~Corwin, J.~Quastel.
\newblock Probability distribution of the free energy of the continuum directed random polymer in $1+1$ dimensions.
\newblock {\em Commun. Pure Appl. Math.}, {\bf 64}:466--537, 2011.

\bibitem{BG}
L.~Bertini, G.~Giacomin.
\newblock  Stochastic Burgers and KPZ equations from particle systems.
\newblock {\em Commun. Math. Phys.}, {\bf 183}:571--607, 1997.

\bibitem{Bethe}
H.~A.~Bethe.
\newblock Zur theorie der metalle. i. eigenwerte und eigenfunktionen der linearen atomkette.
\newblock  {\it Z. Phys.}, {\bf 71}:205, 1931.

\bibitem{BorCor}
A.~Borodin, I.~Corwin.
\newblock Macdonald processes.
\newblock arXiv:1111.4408.


\bibitem{BCS}
A.~Borodin, I.~Corwin, T.~Sasamoto.
\newblock From duality to determinants for q-TASEP and ASEP.
\newblock arXiv:1207.5035.


\bibitem{ICreview}
I.~Corwin.
\newblock The Kardar-Parisi-Zhang equation and universality class.
\newblock arXiv:1106.1596.

\bibitem{Dorlas}
T.~C.~Dorlas.
\newblock Orthogonality and completeness of the Bethe ansatz eigenstates of the nonlinear Schroedinger model.
\newblock  {\it Commun. Math. Physics.}, {\bf 154}:347--376, 1993.

\bibitem{G}
J.~G\"artner.
\newblock Convergence towards Burgers equation and propagation of chaos for weakly asymmetric exclusion process.
\newblock {\em Stoch. Proc. Appl.}, {\bf 27}:233--260, 1988.


\bibitem{GwaSpohn}
L.~H.~Gwa, H.~Spohn.
\newblock Bethe solution for the dynamical-scaling exponent of the noisy Burgers equation.
\newblock {\it Phys. Rev. A}, {\bf 46}:844--854, 1992.

\bibitem{Hahn}
W.~Hahn.
\newblock Beitr\"{a}ge zur Theorie der Heineschen Reihen. Die 24 Integrale der hypergeometrischen q-Differenzengleichung. Das q-Analogon der Laplace-Transformation
\newblock {\it Mathematische Nachrichten}, {\bf 2}:340--379, 1949.


\bibitem{HeckOp}
G.~J.~Heckman, E.~M.~Opdam.
\newblock Yang's system of particles and Hecke algebras.
\newblock {\it Ann. Math.}, {\bf 145}:139--173, 1997.

\bibitem{IS}
T.~Imamura, T.~Sasamoto.
\newblock Current moments of 1D ASEP by duality.
\newblock {\em J. Stat. Phys.},  {\bf 142}:919--930, 2011.

\bibitem{KJ}
K.~Johansson.
\newblock Shape fluctuations and random matrices.
\newblock {\em Commun. Math. Phys.}, {\bf 209}:437--476, 2000.

\bibitem{K}
M.~Kardar.
\newblock Replica-Bethe Ansatz studies of two-dimensional interfaces with quenched random impurities.
\newblock {\it Nucl. Phys. B}, {\bf 290}:582--602, 1987.

\bibitem{KPZ}
K.~Kardar, G.~Parisi, Y.Z.~Zhang.
\newblock  Dynamic scaling of growing interfaces.
\newblock {\em Phys. Rev. Lett.},  {\bf 56}:889--892, 1986.


\bibitem{LL}
E.H.~Lieb, W.~Liniger.
\newblock Exact Analysis of an Interacting Bose Gas. I. The General Solution and the Ground State.
\newblock {\it Phys. Rev. Lett.}, {\bf 130}:1605--1616, 1963.

\bibitem{Lig}
T.~Liggett.
\newblock {\it Interacting particle systems}.
\newblock Spinger-Verlag, Berlin, 2005.


\bibitem{MGP}
J.~MacDonald, J.~Gibbs, A.~Pipkin.
\newblock Kinetics of biopolymerization on nucleic acid templates.
\newblock {\it Biopolymers}, {\bf 6}, 1968.


\bibitem{Oxford}
S.~Oxford.
\newblock {\it The Hamiltonian of the quantized nonlinear Schr\"{o}dinger equation}.
\newblock Ph.D. thesis, UCLA, 1979.

\bibitem{ProSpoComp}
S.~Prolhac, H.~Spohn.
\newblock The propagator of the attractive delta-Bose gas in one dimension.
\newblock {\it J. Math. Phys.}, {\bf 52}:122106, 2011.


\bibitem{SS1}
T.~Sasamoto, H.~Spohn.
\newblock The crossover regime for the weakly asymmetric simple exclusion process.
\newblock {\em J. Stat. Phys.}, {\bf 140}:209--231,  2010.

\bibitem{SchutzMastereq}
G.~M.~Sch\"{u}tz.
\newblock Exact solution of the master equation for the asymmetric exclusion process.
\newblock {\it J. Stat. Phys.}, {\bf 88}:427--445, 1997.

\bibitem{SchutzDuality}
G.~M.~Sch\"{u}tz.
\newblock Duality relations for asymmetric exclusion processes.
\newblock {\em J. Stat. Phys.}, {\bf 86}:1265--1287, 1997.

\bibitem{Spitzer}
F.~Spitzer.
\newblock Interaction of Markov processes.
\newblock {\em Adv. Math.}, {\bf 5}:246--290, 1970.

\bibitem{Sutherland}
B.~Sutherland.
\newblock {\it Beautiful models: 70 years of exactly solvable quantum many-body problems}.
\newblock World Scientific, 2004.

\bibitem{TW1}
C.~Tracy, H.~Widom.
\newblock Integral formulas for the asymmetric simple exclusion process.
\newblock {\em Commun. Math. Phys.}, {\bf 279}:815--844, 2008.
\newblock Erratum: {\em Commun. Math. Phys.} {\bf 304}:875--878, 2011.

\bibitem{TW2}
C.~Tracy, H.~Widom.
\newblock A Fredholm determinant representation in ASEP.
\newblock {\em J. Stat. Phys.}, {\bf 132}:291--300, 2008.

\bibitem{TW3}
C.~Tracy, H.~Widom.
\newblock Asymptotics in ASEP with step initial condition.
\newblock {\em Commun. Math. Phys.}, {\bf 290}:129--154, 2009.


\bibitem{TWBose}
C.~Tracy, H.~Widom.
\newblock The dynamics of the one-dimensional delta-function Bose gas.
\newblock {\it J. Phys. A}, {\bf 41}:485204, 2008.

\bibitem{YangYang}
C.~N.~Yang, C.~P.~Yang.
\newblock One-dimensional chain of anisotropic spin-spin interactions. I. Proof of Bethe's hypothesis for the ground state in a finite system.
\newblock {\it Phys. Rev.} {\bf 150}:321--327, 1966.

\end{thebibliography}
\end{document}